\documentclass[a4paper, 11pt]{article}
\usepackage{amsmath} \usepackage{amsfonts} \usepackage{amssymb}\usepackage{latexsym}
\usepackage[english]{babel}
\usepackage{amscd}
\usepackage[dvips,final]{graphics}
\usepackage{fancyhdr}
\usepackage{locengli,math}
\usepackage{graphicx,pst-plot,psfig,pstricks,pst-node,pst-text,pst-tree,pst-char,pst-coil}

\newcommand{\per}{\textrm{per}}
\newcommand{\F}{\mathcal{F}}
\newcommand{\G}{\mathcal{G}}

\begin{document}

\begin{center}
\textbf{\LARGE{\textsf{Tiling the $(n^2,1)$-De Bruijn graph with $n$ coassociative coalgebras}}} 

\vskip1cm
\parbox[t]{14cm}{\large{
Philippe {\sc Leroux}}\\
\vskip4mm
{\footnotesize
\baselineskip=5mm
Institut de Recherche
Math\'ematique, Universit\'e de Rennes I and CNRS UMR 6625\\
Campus de Beaulieu, 35042 Rennes Cedex, France, philippe.leroux@univ-rennes1.fr}}
\end{center}

\vskip1cm
{\small
\vskip1cm
\baselineskip=5mm
\noindent
{\bf Abstract:}
We construct, via usual graph theory
a class of associative dialgebras as well as a class of coassociative 
$L$-coalgebras, the two classes being related
by a tool from graph theory called the line-extension. 
As a corollary, a tiling of the $n^2$-De Bruijn graph with $n$ (geometric supports of) coassociative coalgebras is obtained. We get, via the tilling of the $(3,1)$-De Bruijn graph, an example of cubical trialgebra defined by J-L Loday and M. Ronco. Other examples are obtained by letting $M_n(k)$ acts on the axioms defining such tilings. Examples of associative products which split into several associative ones are also given.
\section{Introduction}
The field $k$ stands for the real field or the complex field.
Moreover, all the vector spaces considered in this paper will have a finite or
a denumerable basis.
Let us recall
the general setting of this article by summarising
the main steps of our previous works \cite{Coa,dipt,perorb1,tresse}.

\begin{defi}{[Directed graph]}
A {\it{directed graph}} $G$ is a quadruple \cite{Rosen}, $(G_{0},G_{1},s,t)$
where $G_{0}$ and $G_1$ are two denumerable sets respectively called the {\it{vertex set}} and the {\it{arrow set}}.
The two mappings, $s, \ t: G_1 \xrightarrow{} G_0$ are respectively called  {\it{source}} and {\it{terminus}}.
A vertex $v \in G_0$ is a {\it{source}} (resp. a {\it{sink}}) if $t^{-1}(\{v \})$ (resp. $s^{-1}(\{v\})$)
is empty. A graph $G$ is said {\it{locally finite}}, (resp. {\it{row-finite}}) if 
$t^{-1}(\{v\})$ is finite (resp. $s^{-1}(\{v\})$ is finite).
Let us fix a vertex $v \in G_0$.
Define the set $F_{v} :=\{a \in G_{1}, \ s(a)=v \}$. A {\it{weight}} associated with
the vertex $v$ is a mapping $w_v: F_{v} \xrightarrow{} k$.
A directed graph equipped with a family of weights $w := (w_v)_{v \in G_0}$ 
is called a {\it{weighted directed graph}}.
\end{defi}
In the sequel, directed graphs will be supposed locally finite and row finite.
Let us introduce particular coalgebras named $L$-coalgebras \footnote{This notion has been introduced in
\cite{Coa} and developed in \cite{Coa,dipt,perorb1,tresse}.} and explain why this notion is interesting.
\begin{defi}{[$L$-coalgebra]}
A {\it{$L$-coalgebra}} $(L, \Delta, \tilde{\Delta})$ 
over a field $k$ is a $k$-vector space composed of a right part $(L, \Delta)$, where 
$\Delta: L \xrightarrow{} L^{\otimes 2}$, is called the right coproduct and a left part $(L, \tilde{\Delta})$, where $\tilde{\Delta}: L \xrightarrow{} L^{\otimes 2}$, is called the left coproduct such that the
coassociativity breaking equation, 
$(\tilde{\Delta} \otimes id)\Delta = (id \otimes \Delta)\tilde{\Delta}$, is verified.
If $\Delta = \tilde{\Delta}$, the $L$-coalgebra is said  {\it{degenerate}}. A $L$-coalgebra may have two counits, the right counit  
$\epsilon: L \xrightarrow{} k$, verifying $ (id \otimes \epsilon)\Delta = id$
and the left counit $\tilde{\epsilon}: L \xrightarrow{} k$, verifying $ ( \tilde{\epsilon} \otimes id)\tilde{\Delta} = id.$ A $L$-coalgebra is
said {\it{coassociative}} if its two coproducts are coassociative. In this case the equation, $(\tilde{\Delta} \otimes id)\Delta = (id \otimes \Delta)\tilde{\Delta}$, is called the {\bf{entanglement equation}} and we will say that its right part $(L, \Delta)$ is entangled to its left part $(L, \tilde{\Delta})$.
Denote by $\tau$, the {\it{transposition}} mapping, i.e. $L^{ \otimes 2} \xrightarrow{\tau} L^{ \otimes 2}$ such that $\tau(x \otimes y) = y \otimes x$, for all $x,y \in 
L$.
The $L$-coalgebra $L$ is said to be {\it{$L$-cocommutative}} if for all $v \in L$, $(\Delta - \tau\tilde{\Delta})v = 0$. 
A \textit{$L$-bialgebra} 
(with counits $\epsilon, \ \tilde{\epsilon}$), is a $L$-coalgebra (with counits) and an unital algebra such that 
its coproducts and counits are homomorphisms.
A \textit{$L$-Hopf algebra}, $H$, is a $L$-bialgebra with counits equipped with right and left antipodes
$S, \ \tilde{S}: H \xrightarrow{} H $, such that:
$
m(id \otimes S)\Delta = \eta \epsilon$ and 
$m(\tilde{S} \otimes id)\tilde{\Delta} = \eta \tilde{\epsilon}$
or
$m(S \otimes id)\Delta = \eta \epsilon$ and
$m(id \otimes \tilde{S})\tilde{\Delta} = \eta \tilde{\epsilon}.$
\end{defi}

Let $G=(G_0, G_1,s,t)$ be a directed graph equipped with a family of weights $(w_v)_{v \in G_0}$. Let us consider the free $k$-vector space $kG_0$ generated by $G_0$.
The set $G_1$ is then viewed as a subset of $(kG_0)^{\otimes 2}$ by identifying $a \in G_1$ with $s(a) \otimes t(a)$. The 
mappings source and terminus are then linear mappings still called 
source and terminus $s,t: \ (kG_0)^{\otimes2} \xrightarrow{} kG_0$, such that $s(u \otimes v) = u$ and $t(u \otimes v) = v$, for all $u,v \in G_0$. The family of weights $(w_v: F_v \xrightarrow{} k)_{v \in G_0}$, is then viewed as a family of linear mappings.
Let $v \in G_0$. Define the right coproduct $\Delta_M: kG_0 \xrightarrow{} (kG_0)^{\otimes 2}$, such that
$\Delta_M(v) := \sum_{i: a_i \in F_v} \ w_v(a_i) \ v \otimes t(a_i)$ and the left coproduct
$\tilde{\Delta}_M: kG_0 \xrightarrow{} (kG_0)^{\otimes 2}$, such that $\tilde{\Delta}_M (v) := \sum_{i: a_i \in P_v} \ w_{s(a_i)}(a_i) \ s(a_i) \otimes v$, where $P_{v} $ is the set $\{a \in G_{1}, \ t(a)=v \}$. 
With these definitions the $k$-vector space $kG_0$ is a $L$-coalgebra called a finite Markov $L$-coalgebra since its coproducts $\Delta_M$ and  $\tilde{\Delta}_M$ verify
the
coassociativity breaking equation 
$(\tilde{\Delta}_M \otimes id)\Delta_M = (id \otimes \Delta_M)\tilde{\Delta}_M$. This particular coalgebra is called in addition finite Markov ($L$-coalgebra)
because for all $v \in G_0$, the sets $F_v$ and $P_v$ are finite and the coproducts have the form $\Delta_M(v) := v \otimes \cdots $ and $\tilde{\Delta}_M(v) := \cdots \otimes v $.

Assume we consider the Markov $L$-coalgebra just described and
associate
with each tensor product $\lambda u \otimes v$, where $\lambda \in k$ and $u,v \in G_0$, appearing in the definition of the coproducts, a directed arrow
$u \xrightarrow{\lambda} v$. The weighted directed graph so obtained, called the {\it{geometric support}} of this $L$-coalgebra, is up to a graph isomorphism \footnote{A {\it{graph isomorphism}} $f: G \xrightarrow{} H$ between two graphs $G$ and $H$ is a pair of bijection $f_0: G_0 \xrightarrow{} H_0$ and $f_1: G_1 \xrightarrow{} H_1$ such that
$f_0(s_G(a))=s_H(f_1(a))$ and $f_0(t_G(a))=t_H(f_1(a))$ for all $a \in G_1$. All the directed graphs in this formalism will be considered up to a graph isomorphism.}, the directed graph we start with. Therefore, general $L$-coalgebras generalise the notion of weighted directed graph. If $(L,\Delta,\tilde{\Delta})$
is a $L$-coalgebra generated as a $k$-vector space by a set $L_0$, then its geometric support 
$Gr(L)$ is a directed graph with vertex set $Gr(L)_0=L_0$ and with arrow set $Gr(L)_1$, the 
set of tensor products $u \otimes v$, with $u,v \in L_0$, appearing
in the definition of the coproducts of $L$. As a coassociative coalgebra is a particular $L$-coalgebra, we naturally construct its directed graph.
We draw attention to the fact that a directed graph can be the
geometric support of different $L$-coalgebras.
\begin{exam}{}
The directed graph:
\begin{center}
\includegraphics*[width=3.5cm]{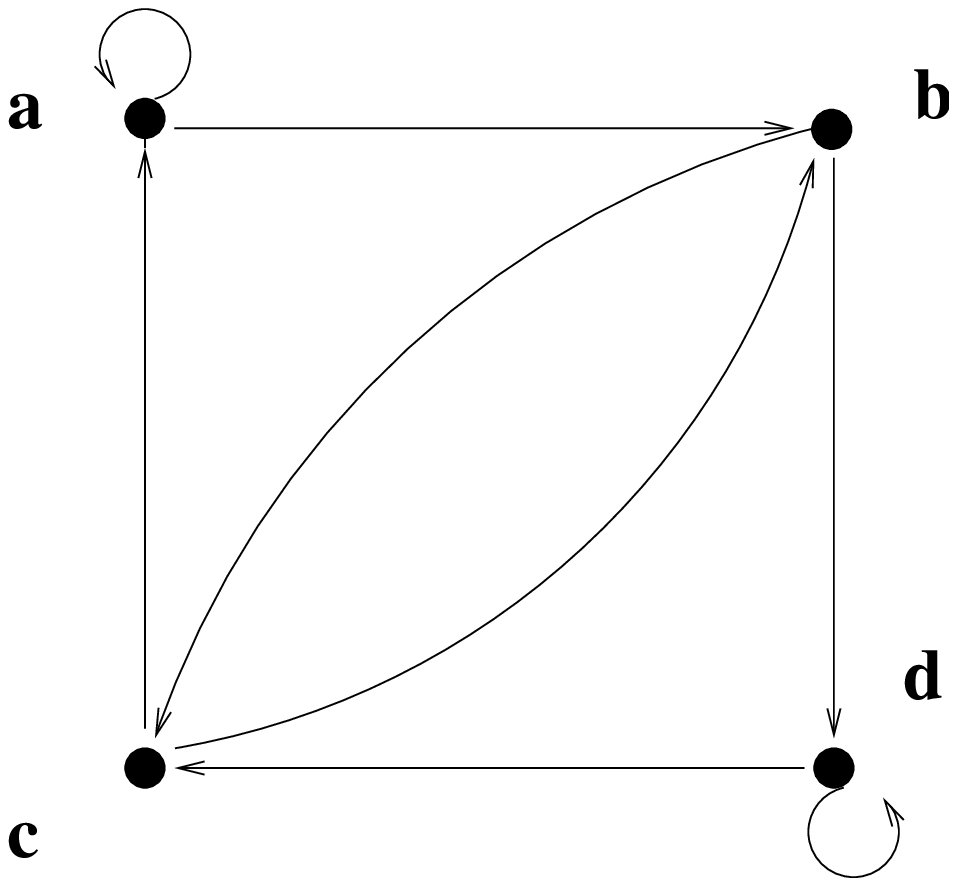}
\end{center}
\end{exam}
is the geometric support associated with the degenerate $L$-coalgebra or coassociative coalgebra $\F$, generated by $a,b,c$ and $d$, as a $k$-vector space, and
described by the following coproduct:
$
\Delta a = a \otimes a + b \otimes c, \ \
\Delta b = a \otimes b + b \otimes d, \ \
\Delta c = d \otimes c + c \otimes a, \ \
\Delta d = d \otimes d + c \otimes b
$
and the geometric support of the finite Markov $L$-coalgebra, generated by $a,b,c$ and $d$, as a $k$-vector space, and described by the right coproduct:
$
\Delta_M a = a \otimes (a + b), \ \
\Delta_M b = b \otimes (c + d), \ \
\Delta_M c = c \otimes  (a +b), \ \
\Delta_M d = d \otimes (c + d)
$
and the left coproduct: 
$
\tilde{\Delta}_M a = (a+c) \otimes a,  \ \
\tilde{\Delta}_M b = (a +c) \otimes b , \ \
\tilde{\Delta}_M c = (b + d) \otimes c , \ \
\tilde{\Delta}_M d = (b+ d) \otimes d.
$
\NB
Let $(\G, \Delta_M, \tilde{\Delta}_M)$ be a finite Markov $L$-coalgebra generated as a $k$-vector space by a set $\G_0$. If the family of weights $(w_v)_{v \in \G_0 }$, $(\tilde{w}_v)_{v \in \G_0 }$
used for describing right and left
coproducts take values into $\mathbb{R}_+$ and
if the right counit $\epsilon: v \mapsto 1$ exists, then the geometric
support associated with $(\G, \Delta_M, \tilde{\Delta}_M)$
is a directed graph equipped with a family of probability vectors described by $(w_v)_{v \in \G_0 }$. 

Before going on, let us interpret what represents the $L$-cocommutativity
in the case of a finite Markov $L$-coalgebra. A directed graph is said bi-directed if for any arrow from a vertex $v_1$ to a vertex $v_2$, there exists an arrow from $v_2$ to $v_1$. To take into account the bi-orientation of a directed graph in an algebraic way, we first embed this directed graph into the finite
Markov $L$-coalgebra described above. We then notice that a directed graph is bi-directed if and only if $\Delta_M = \tau \tilde{\Delta}_M$. Therefore, in this algebraic framework, we are naturally
led to consider the $L$-cocommutator $\ker (\Delta - \tau \tilde{\Delta})$. Dualizing this formula leads to consider an ($L$-)algebra $D$ equipped with two products $\vdash$ and $\dashv$, verifying $(x \vdash y) \dashv z = x \vdash (y \dashv z)$, $x,y,z \in D$, and to consider
the particular commutator $[x,y]:= x \dashv y - y \vdash x$. The bracket, $[ -,z]$, verifies the
analogue of the ``Jacobi identity'', called the Leibniz identity, i.e. 
$[[x,y ] ,z ] = [[ x,z ] ,y ] + [ x,[ y  ,z ]],$
if $D$ is an algebra called an associative dialgebra \cite{Loday}. 

Another motivation concerning associative dialgebras is the following.
In a long-standing project whose ultimate aim is to study periodicity phenomena in algebraic $K$-theory,
J-L. Loday in \cite{Loday}, and J-L. Loday and M. Ronco in  \cite{LodayRonco} introduce several kind of algebras, one of which is 
the ``non-commutative Lie algebras'', called {\it{Leibniz algebras}}. Such algebras $D$ are described by a bracket $[ -,z]$ verifying the so called Leibniz identity:
$$[[x,y ] ,z ] = [[ x,z ] ,y ] + [ x,[ y  ,z ]].$$
When the bracket is skew-symmetric, the Leibniz identity becomes the
Jacobi identity and the Leibniz algebra turns out to be a Lie algebra.
A way to construct such Leibniz algebra is to start with an {\it{associative dialgebra}}, that is a $k$-vector space $D$ equipped with two associative products,
$\vdash$ and $\dashv$, such that for all $x,y,z \in D$ 
\begin{enumerate}
\item{$x \dashv (y \dashv z) = x \dashv (y \vdash z),$}
\item{$(x \vdash y) \dashv z = x \vdash (y \dashv z),$}
\item{$(x \dashv y) \vdash z = (x \vdash y) \vdash z.$}
\end{enumerate}
The associative dialgebra is then a Leibniz algebra by defining the bracket
$[x,y ] := x \dashv y - y \vdash x$, for all $x,y \in D$. 
The operad associated with associative dialgebras is then Koszul dual to the operad associated with dendriform algebras,
a {\it{dendriform algebra}} $Z$ being a $k$-vector space equipped with
two binary operations,
$ \prec \ , \ \succ: \ Z \otimes Z \xrightarrow{} Z,$
satisfying the following axioms:
\begin{enumerate}
\item {$(a \prec b) \prec c = a \prec (b \prec c) + a \prec (b \succ c),$ }
\item {$(a \succ b) \prec c = a \succ (b \prec c),$ }
\item {$(a \prec b) \succ c + (a \succ b) \succ c = a \succ (b \succ c).$ }
\end{enumerate}
This notion dichotomizes the notion of associativity since the product
$a *b = a \prec b + a \succ b$, for all $a,b \in Z$ is associative. Otherwise stated, the associative product $*$ splits into two  operations $\prec$ and $\succ$.

The first important result of this paper is the construction, via Markov $L$-coalgebras, i.e. via usual graph theory
of a class of $L$-cocommutative and coassociative codialgebras as well as a class of coassociative 
$L$-coalgebras, the two classes being related
by a tool from graph theory called the line-extension. 
The second one is the construction of a tiling of the $(n^2,1)$-De Bruijn graph with $n$ (geometric supports of) coassociative coalgebras. As a corollary, we obtain examples of cubical trialgebras, a notion developed by J-L Loday and M. Ronco in \cite{LodayRonco} and splittings of associative products into several associative ones.

Let us briefly introduce the organisation of the paper.
In section 2, we display the notion of coassociative co-dialgebras and recall the definition of the De-Bruijn graphs. We prove that the $(2,1)$-De Bruijn graph,
viewed as a Markov co-dialgebra, yields, by line-extension, the geometric support of the coassociative coalgebra $\F$, generated by say $a,b,c$ and $d$, as a $k$-vector space, and
described by the following coproduct:
$
\Delta a = a \otimes a + b \otimes c, \ \
\Delta b = a \otimes b + b \otimes d, \ \
\Delta c = d \otimes c + c \otimes a, \ \
\Delta d = d \otimes d + c \otimes b
$.
Inspired by a previous work \cite{perorb1}, relationships between the Markovian coproduct of the $(2,1)$-De Bruijn graph and the coassociative coproduct
of $\F$ are also given. In section 3, we prove the existence of another structure associated with $\F$. This structure is also a coassociative coalgebra,
called by convention the left structure of $\F$.
As the two structures of $\F$ are entangled by the entanglement equation , we conclude that $\F$ is a coassociative $L$-coalgebra.
Moreover, gluing their two associated geometric supports yields the $(4,1)$-De Bruijn graph. Since the intersection of the two arrow sets of their geometric supports is
empty, we assert that the $(4,1)$-De Bruijn graph can be tiled by the two geometric supports of the coassociative $L$-coalgebra $\F$.

The consequences of such a left structure on $Sl_q(2)$ are also explored.
If $\F$ stands for $Sl_q(2)$, we notice that the usual algebraic relations of the Hopf algebra $Sl_q(2)$ do not embed this new structure
into a bialgebra. However, it allows us to construct a map which verify the same axioms as an antipode map. We yield also
a left structure for $SU_q(2)$.

This work ends, in section 4,
by showing that the $(n^2,1)$-De Bruijn graph can also be tiled with $n$ (geometric supports of
) coassociative coalgebras. As a consequence, we give examples of cubical trialgebras and put forward an important notion which is the achirality of a $L$-coalgebra.
we give also
Leibniz co-derivative on these directed graphs.
\section{On the De Bruijn graph families}
\label{bruijn}
The description of weighted directed graphs and coassociative coalgebras can be embedded into
the $L$-coalgebra framework. In \cite{Coa}, we establish an important theorem. 
\begin{prop}
Any coassociative coalgebra $(C, \Delta_C)$, with a group-like element can be embedded into a non degenerate $L$-coalgebra.
\end{prop}
\Proof
Let $(C, \Delta_C)$ be a coassociative coalgebra. Suppose $e$ is a group-like element, i.e. $\Delta_C e = e \otimes e$.
Define the coassociative coproducts $\delta(c):=c \otimes e$ and $\tilde{\delta}(c):=e \otimes c$ for all $c \in C$. 
Define also the linear map
$\overrightarrow{d}: \ C \xrightarrow{} C \otimes C \ \textrm{such that} \ \overrightarrow{d}(c)=\Delta c - \delta_f(c)$
and the linear map
$\overleftarrow{d}: \ C \xrightarrow{} C \otimes C \ \textrm{such that} \ \overleftarrow{d}(c)=  \Delta c -\tilde{\delta}_f(c).$
These linear maps, $\overleftarrow{d}$ and $\overrightarrow{d}$, turn the coassociative coalgebra $(C,\Delta)$ into a non degenerate $L$-coalgebra $(C,\overrightarrow{d},\overleftarrow{d})$.
\eproof 

\noindent
Let $C$ be a bialgebra with unit $e$. The two
new coproducts, $\overleftarrow{d}, \ \overrightarrow{d}: C \xrightarrow{} C^{\otimes 2}$, turn out to be Ito derivative \footnote{Let $A$ be an associative algebra with unit $e$ and $f: A \xrightarrow{} A$ be a linear map. The map $f$ is said to be an Ito derivative if $f(e) =0$ and 
$f(xy) = f(x)f(y) + xf(y) + f(x)y$, for all $x,y \in A$.}.
Starting with a Markov $L$-coalgebra $\G$, we can
recover a similar theorem.
However, whereas in the coassociative coalgebra case, these new two coproducts map ``the vertex set'' $C$ into the ``arrow sets'' $C^{\otimes 2}$,
the new coproducts, $\overleftarrow{d}_\G$ and $\overrightarrow{d}_\G$, of a Markov $L$-coalgebra $\G$ map ``the arrow set'' $C^{\otimes 2}$ into the `` paths of lenght 2'' modelised by $C^{\otimes 3}$.
This observation suggests the necessity of studying the line-extension of directed graphs, a notion defined in the sequel.
\subsection{Line-extension of directed graphs}
The aim of this subsection is to show that some De Bruijn directed graphs, seen as Markov $L$-coalgebras are also coassociative co-dialgebras and that
the line-extension of these directed graphs can be viewed as geometric supports associated with some well-known coassociative coalgebras.
We start with two definitions.
\begin{defi}{[De Bruijn graph]}
A $(p,n)$ De-Bruijn sequence on the alphabet $\Sigma = \{ a_1, \ldots, a_p \}$ is a sequence $(s_1, \ldots, s_{m})$ of $m = p^n$ elements
$s_i \in \Sigma$ such that subsequences of length $n$ of the form $(s_i, \ldots, s_{i+n -1})$  are distinct,
the addition of subscripts being done modulo $m$.
A {\it{$(p,n)$-De Bruijn graph}} is a directed graph
whose vertices correspond to all possible strings $s_1s_2 \ldots s_n$ of $n$ symbols from $\Sigma$.
There are $p$ arcs leaving the vertex $s_1s_2 \ldots s_n$ and leading to the adjacent node $s_2s_3 \ldots s_n \alpha$, $\alpha \in \Sigma$.
Therefore the {\it{$(p,1)$-De Bruijn graph}} is the directed graph with $p$ vertices, complete, with a loop at each vertex.

\noindent
Fix $n \geq 1$.
Let $D_{(n,1)}=((D_{(n,1)})_0, (D_{(n,1)})_1, s,t)$ be the $(n,1)$-De Bruijn graph. Its natural Markov $L$-coalgebra is defined as follows.
Denote by $v_i$, $1 \leq i \leq n$, the vertices of 
$D_{(n,1)}$. Embed $(D_{(n,1)})_0$ into its free $k$-vector space. Define the coproducts $\Delta_M , \ \tilde{\Delta}_M: k(D_{(n,1)})_0 \xrightarrow{}k(D_{(n,1)})_0^{\otimes 2}$, such that for all $i$, $\Delta_M v_i:= v_i \otimes \sum_j v_j$ and 
$\tilde{\Delta}_M v_i:= \sum_j v_j \otimes  v_i$.
There are obvious left and right counits,
$ \tilde{\epsilon}(v_i)=\epsilon(v_i)=\frac{1}{n}$, for all $v_i \in (D_{(n,1)})_0$.
\end{defi}
\begin{defi}{[Line-extension]}
The {\it{line-extension}} of a directed graph $G$,
with vertex set $G_0= \{j_1, \ldots, j_n \}$ and arrow set $G_1 \subseteq G_0 \times G_0$, denoted by $E(G)$,
is the directed graph with vertex set $E(G)_0 := G_1$ and the arrow set $E(G)_1 \subseteq G_1 \times G_1$ defined by
$((j_k, j_l), (j_e, j_f))$ iff $j_l = j_e$. 
\end{defi}
The notion of associative dialgebra, due to J-L. Loday \cite{Loday}, is a notion which generalizes the notion of algebra. Associative dialgebras,
via the notion of dendriform algebra \footnote{The operads Dias (associative dialgebra) and
Dend (dendriform algebra) are dual in the operad sense, see \cite{Loday}.}, are closely related to (planar binary) trees,
which became an important tool in quantum field theory \cite{Kreimer,CoKreimer}.
Here, we are interested in the notion of coassociative co-dialgebra.
\begin{defi}{[Coassociative co-dialgebra of degree $n$]}
Motivated by line-extension of (geometric support of) Markov $L$-coalgebra, we define in \cite{Coa}, (Markov) $L$-coalgebra of degree $n$, $n>0$.
Similarly,
let $\Delta_n$ and $\tilde{\Delta}_n$ be two $n$-linear mappings $D^{\otimes n} \xrightarrow{}D^{\otimes n+1}$, where $D$ is a $k$-vector space. 
The $k$-vector space $(D,\Delta_n, \tilde{\Delta}_n)$ is said to be a {\it{coassociative co-dialgebra}} of degree $n$ if the following axioms are verified:
\begin{enumerate}
\item {$\Delta_n$ and $\tilde{\Delta}_n$ are coassociative,}
\item{$(id \otimes \Delta_n)\Delta_n = (id \otimes \tilde{\Delta}_n) \Delta_n $,}
\item{$(\tilde{\Delta}_n \otimes id) \tilde{\Delta}_n = (\Delta_n  \otimes id) \tilde{\Delta}_n$,}
\item{$(\tilde{\Delta}_n \otimes id) \Delta_n = ( id \otimes \Delta_n) \tilde{\Delta}_n$.}
\end{enumerate}
As in the $L$-coalgebra case, such a space may have a right counit
$\epsilon_n: D^{\otimes n} \xrightarrow{} D^{\otimes n-1} \  \textrm{such that:} \  (id \otimes \epsilon_n)\Delta_n = id$
and a left counit
$\tilde{\epsilon}_n: D^{\otimes n} \xrightarrow{} D^{\otimes n-1} \  \textrm{such that:} \  ( \tilde{\epsilon}_n \otimes id)\tilde{\Delta}_n = id.$
By convention $D^{\otimes 0} := k$.
A coassociative co-dialgebra of degree 1 will be also called a coassociative co-dialgebra.
\end{defi}
\begin{prop}
A $k$-vector space $C$ equipped with two coproducts $\delta_f$ and $\tilde{\delta}_f$ such that $\delta_f(c) = c \otimes e$ and $\tilde{\delta}_f(c)= e \otimes c$, for all $c \in C$ is a Markov coassociative codialgebra.
\end{prop}
\Proof
Straightforward since $\delta_f(e) = \tilde{\delta}_f(e)$.
\eproof
\begin{exam}{[The flower graph]}
A unital algebra, $A$ with unit 1, carries a non trivial Markov $L$-bialgebra, obtained from the equality, $(1 \cdot a)\cdot 1 =1 \cdot( a\cdot 1)$.
The coproducts are $\delta_f(a) = a \otimes 1$ and $\tilde{\delta}_f(a)= 1 \otimes a$, $a \in A$. Its geometric support is called the flower graph.
\begin{center}
\includegraphics*[width=4cm]{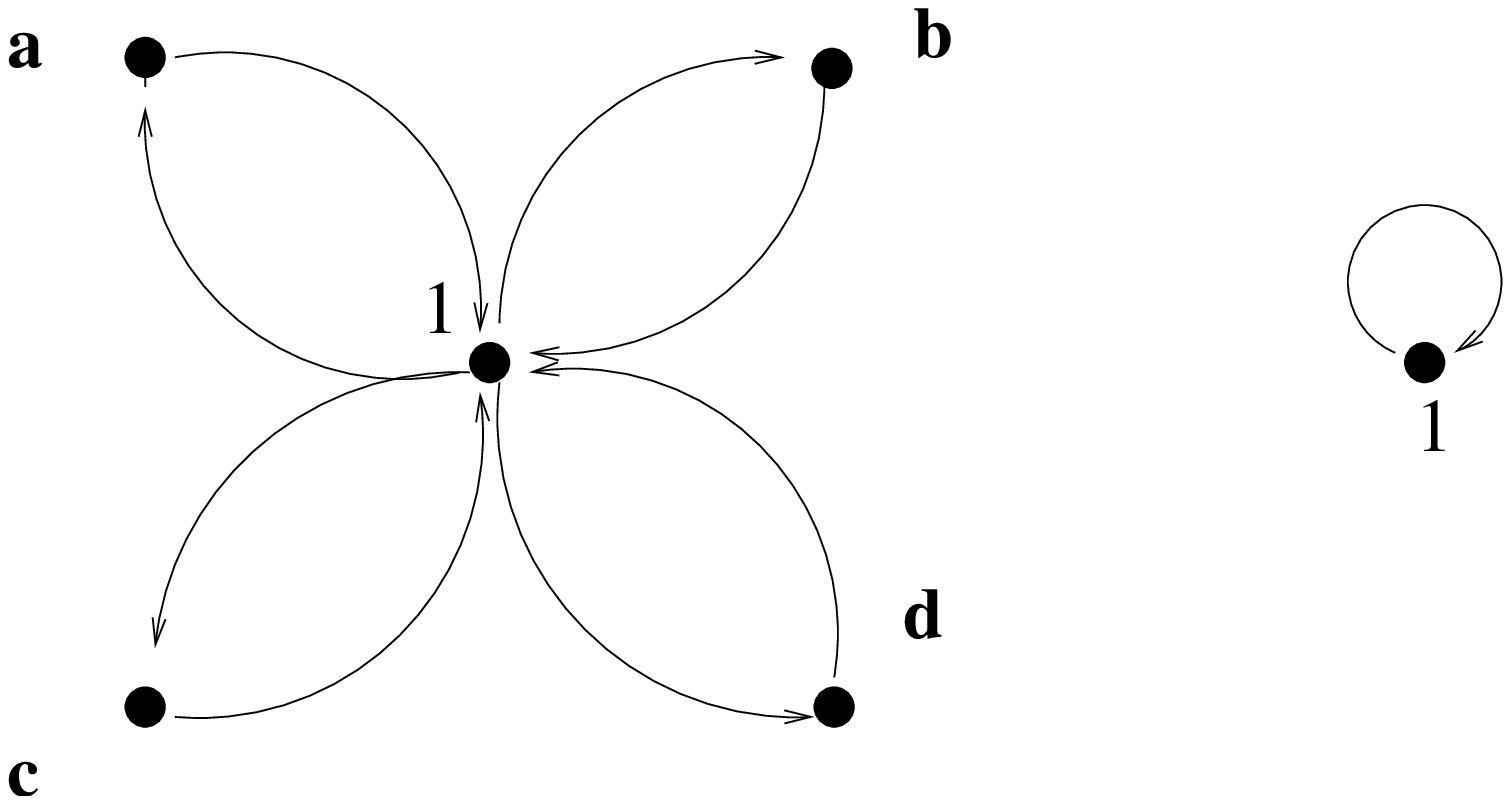}

\begin{scriptsize} \textbf{Example of geometric support associated with an algebra $k \bra a,b,c,d \ket \oplus k1 $.} \end{scriptsize}
\end{center}
\end{exam}
\begin{theo}
\label{aze}
Any coassociative coalgebra $(C, \Delta)$, (respectively bialgebra, Hopf algebra), can be embedded into a $L$-coalgebra of degree $n$, (respectively, $L$-bialgebra of degree $n$,
$L$-Hopf algebra \footnote{Set $id_{n}:= \underbrace{id \otimes id \otimes \ldots \otimes id}_{n} $, $n >0$. A $L$-Hopf algebra of degree $n$, $(H,\Delta_H,\tilde{\Delta}_H)$, is by definition a $L$-bialgebra of degree $n$, equipped with right and left counits, $\tilde{\epsilon}_H$, $\epsilon_H$ of degree $n$,
such that its antipodes $S, \tilde{S}: \ H \xrightarrow{} H $ verify $(id_{n-1} \otimes m) (id_{n} \otimes S) \Delta_H = \eta_n	 \epsilon_H$ and $
(m \otimes id) (\tilde{S} \otimes id_{n}) \tilde{\Delta}_H = \tilde{\eta}_n \tilde{\epsilon}_H,$
with $\eta_n, \tilde{\eta}_n: H^{\otimes (n-1)} \xrightarrow{} H^{\otimes n}$ such that $\eta_n(h) := h \otimes 1_H$ and $\tilde{\eta}_n(h) := 1_H \otimes h$, $h \in H^{\otimes (n-1)}$.} of degree $n$), with $n>1$.
\end{theo}
\Proof
Fix $n>1$.
Consider $\Delta$, the coproduct of such a coassociative coalgebra $C$. Set:
$id_{n}:= \underbrace{id \otimes id \otimes \ldots \otimes id}_{n} $ and $\Delta_n := id_{n-1} \otimes \Delta$
and $\tilde{\Delta}_n :=\Delta \otimes id_{n-1}$. The right and left coassociative coproducts
$\Delta_n $ and $\tilde{\Delta}_n$ map $C^{\otimes n}$ into $C^{\otimes n+1}$ and the entanglement equation is realised. The right and left counits are $\epsilon_n := id_{n-1} \otimes \epsilon$ and
$\tilde{\epsilon}_n:= \epsilon \otimes id_{n-1}$, since for instance $(id \otimes \epsilon_n)\Delta_n = id_{n}$. If the counit and the coproduct of $C$ are unital homomorphisms, so are
the new coproducts and counits. If $C$ is a Hopf algebra with antipode $s$, we embed $C$ into a $L$-Hopf algebra of degree $n$ since  $(id_{n-1} \otimes m)(id_{n} \otimes s)\Delta_n= \eta_n \epsilon_n$ and $(m \otimes id_{n-1})(s \otimes id_{n} )\Delta_n= \tilde{\eta}_n \tilde{\epsilon}_n$.
\eproof
\begin{theo}
If $n=2$, any coassociative coalgebra $(C,\Delta)$, can be viewed as a coassociative co-dialgebra of degree 2.
\end{theo}
\Proof
Straightforward by using the definition of $\Delta_2$ and the fact that $\Delta$ is coassociative.
\eproof
\begin{prop}
Let $D_{(n,1)}=((D_{(n,1)})_0, \ (D_{(n,1)})_1)$, be the $(n,1)$-De Bruijn graph and consider its free $k$-vector space $k(D_{(n,1)})_0$.
The Markovian coproducts of the Markov $L$-coalgebra $k(D_{(n,1)})_0$ associated with the $(n,1)$-De Bruijn graph define a coassociative co-dialgebra.
\end{prop}
\Proof
Straightforward.
\eproof
\begin{prop}
\label{tyu}
There exists a coassociative coalgebra whose geometric support is the line-extension of the $(n,1)$-De Bruijn graph.
\end{prop}
\Proof
Let us denote the arrow emerging from a given vertex $v_i$ to a vertex $v_j$, with $i,j = 1, \ldots, n$ of the $(n,1)$-De Bruijn graph $D_{(n,1)}$ by $a_{ij}$.
The new vertices of $E(D_{(n,1)})$ are denoted by $a_{ij}$ and the arrows are denoted by $((il), (lj))$. Consider the free $k$-vector space generated by the set $E(D_{(n,1)})_0:=\{a_{ij}; \ i,j =1 \ldots n \}$ and define
$\Delta a_{ij}= \sum_l \ a_{il}\otimes a_{lj}$, this coproduct is coassociative and the geometric support associating with
the coassociative coalgebra $(kE(D_{(n,1)})_0, \Delta)$ is easily seen to be $E(D_{(n,1)})$. It has an obvious counit, $a_{ij} \mapsto 0$ if
$i \not= j$ and $a_{ij} \mapsto 1$ otherwise.
\eproof
\begin{coro}
Recall that $\F$ is generated by $a,b,c$ and $d$, as a $k$-vector space, and
described by the following coproduct:
$
\Delta a = a \otimes a + b \otimes c, \ \
\Delta b = a \otimes b + b \otimes d, \ \
\Delta c = d \otimes c + c \otimes a, \ \
\Delta d = d \otimes d + c \otimes b.
$
The line-extension of the $(2,1)$-De Bruijn graph can be equipped with the coassociative coproduct associated with
the coalgebra $\F$.
\end{coro}
\begin{prop}
Let $B$ be a (Markov) coassociative co-dialgebra with coproducts $\tilde{\Delta}$ and $\Delta$ and $C$ be a coassociative coalgebra
with coproduct $\Delta_C$. Then $B \otimes C$
is a coassociative co-dialgebra with coproduct $\delta_{B \otimes C} := (id \otimes \tau \otimes id) \Delta \otimes \Delta_C $ and
$\tilde{\delta}_{B \otimes C} := (id \otimes \tau \otimes id) \tilde{\Delta} \otimes \Delta_C $.
Similarly, $C \otimes B$ is a coassociative co-dialgebra with  $\delta_{C \otimes B} := (id \otimes \tau \otimes id) \Delta_C \otimes \Delta$ and
$\tilde{\delta}_{C \otimes B} := (id \otimes \tau \otimes id)    \Delta_C  \otimes \tilde{\Delta}$.
\end{prop}
\Proof
Straightforward.
\eproof

\noindent
We end this section on De-Bruijn graphs and coassociative codialgebras, by constructing another family of coassociative co-dialgebra. The idea is to give the $(n,1)$-De Bruijn graph an attractor r\^ole.
\begin{prop}
Let $n$ and $m$ be two integers different from zero. Let $D$ be the $k$-vector space $k\bra x_1, \ldots, x_m \ket \oplus k\bra \alpha_1, \ldots, \alpha_n \ket $. Define for all $i = 1, \ldots, m$, $\Delta x_i = \sum_{j=1}^n x_i  \otimes \alpha_j$,
$\tilde{\Delta} x_i = \sum_{j=1}^n  \alpha_j  \otimes x_i$ and
for all $j = 1, \ldots, n$, $\Delta  \alpha_j = \sum_{p=1}^n   \alpha_j \otimes \alpha_p$, $\tilde{\Delta} \alpha_j  = \sum_{p=1}^n  \alpha_p  \otimes  \alpha_j $.
These Markovian coproducts embed $D$ into a coassociative co-dialgebra.
\end{prop}
\Proof
Straightforward.
\eproof
\Rk
Consider the $k$-vector space of linear maps $L(D, A)$ which map $(D,\Delta,\tilde{\Delta})$, a coassociative coalgebra into an associative algebra $A$ equipped
with a product $m$. The space $L(D, A)$ is then embedded into an associative dialgebra by defining for all $f, g \in L(D, A)$, the two convolution products:
$f \dashv g := m(f \otimes g) \Delta$ and $ f \vdash g := m(f \otimes g) \tilde{\Delta}$. 
\subsection{Relationships between the $(2,1)$-De Bruijn $L$-coalgebra and $\F$}
Motivated by a previous work \cite{perorb1}, 
the aim of this subsection is to study the relationships between the $(2,1)$-De Bruijn graph, seen as a Markov $L$-coalgebra and its
line-extension seen as the (geometric support of the) coassociative coalgebra $\F$.
Let $X,Y$ be two non commuting operators, we consider the non commutative algebra $A = k \bra X, \ Y \ket \oplus k1$, with $XY = \eta \ YX$ and $\eta \in k \setminus \{ 0 \}$.
We equip $A$ with the following Markovian coproducts
$\Delta_M X = X \otimes X + X \otimes Y$ and $\Delta_M Y = Y \otimes X + Y \otimes Y$,
$\tilde{\Delta}_M X =Y \otimes X + X \otimes X$, $\tilde{\Delta}_M Y = X \otimes Y + Y \otimes Y$ and set
$ \tilde{\Delta}_M 1 = \Delta_M 1 = 1 \otimes 1$.
Set $a = X \otimes X, \ b= X \otimes Y, \ c= Y \otimes X, d = Y \otimes Y$. We would like to find operators which
give the coassociative coproduct of $\F$ from the coassociative co-dialgebra $A$. For that, we define:
$$ \overrightarrow{\Delta}_M: A \xrightarrow{} A^{\otimes 2} \times A^{\otimes 2}, \ \
X \mapsto  \overrightarrow{\Delta}_M(X) = \begin{pmatrix}
 X \otimes X \\
X \otimes Y
\end{pmatrix},
Y \mapsto  \overrightarrow{\Delta}_M(Y) = \begin{pmatrix}
 Y \otimes Y \\
Y \otimes X
\end{pmatrix},
1 \mapsto  \overrightarrow{\Delta}_M(1) = \begin{pmatrix}
 1 \otimes 1 \\
1 \otimes 1
\end{pmatrix},
$$
and
$$ \overrightarrow{\tilde{\Delta}}_M: A \xrightarrow{} A^{\otimes 2} \times A^{\otimes 2}, \ \
X \mapsto  \overrightarrow{\tilde{\Delta}}_M(X) = \begin{pmatrix}
 X \otimes X \\
Y \otimes X
\end{pmatrix},
Y \mapsto  \overrightarrow{\tilde{\Delta}}_M(Y) = \begin{pmatrix}
 Y \otimes Y \\
X \otimes Y
\end{pmatrix},
1 \mapsto  \overrightarrow{\tilde{\Delta}}_M(1) = \begin{pmatrix}
 1 \otimes 1 \\
1 \otimes 1
\end{pmatrix}.
$$
If we
define the bilinear map $ \bra \cdot \ ; \ \cdot \ket_*: \ (A^{\otimes 2} \times A^{\otimes 2})^{\times 2} \xrightarrow{} A^{\otimes 2},
((z_1, z_2), (z_3, z_4)) \mapsto z_1z_3 + z_2z_4$, we recover
$\Delta_M(X) = \bra \overrightarrow{\Delta}_M(X), \overrightarrow{\Delta}_M(1) \ket_*$ and
$\Delta_M(Y) = \bra \overrightarrow{\Delta}_M(Y), \overrightarrow{\Delta}_M(1) \ket_*$.
Define the bilinear map $\square: \ A^{\otimes 2} \times A^{\otimes 2} \xrightarrow{} A^{\otimes 2} \otimes A^{\otimes 2}$ such that:
$$(y_1 \otimes y_2) \square (y_3 \otimes y_4) := (id \otimes \tau \otimes id) ((y_1 \otimes y_2) \otimes (y_3 \otimes y_4)
:=(y_1 \otimes y_3) \otimes (y_2 \otimes y_4). $$ Define also $ \bra \cdot \ ; \ \cdot \ket$ and $\per: (A^{\otimes 2} \times A^{\otimes 2})^{\times 2} \xrightarrow{} A^{\otimes 2} \otimes A^{\otimes 2}$
such that $((z_1, z_2), (z_3, z_4)) \mapsto z_1 \square z_3 + z_2 \square z_4$
and the ``permanent'' $\per((z_1, z_2), (z_3, z_4)) := z_1 \square z_4 + z_2 \square z_3$.
Let us express the relations between Markovian coproducts of the $(2,1)$-De Bruijn graph and the coassociative coproduct of $\F$.
\begin{prop}
\label{markov} 
The relations are: 
$\bra \overrightarrow{\Delta}_M(X), \overrightarrow{\tilde{\Delta}}_M(X) \ket = \Delta(a)$,
$\bra \overrightarrow{\Delta}_M(Y), \ \overrightarrow{\tilde{\Delta}}_M(Y) \ket = \Delta(d)$, \\
$\per(\overrightarrow{\Delta}_M(X), \overrightarrow{\tilde{\Delta}}_M(Y))=\Delta (b)$,
$\per(\overrightarrow{\Delta}_M(Y), \overrightarrow{\tilde{\Delta}}_M(X))=\Delta (c).$
\end{prop}
\Proof
Straightforward. For instance, $\bra \overrightarrow{\Delta}_M(X), \overrightarrow{\tilde{\Delta}}_M(X) \ket = (X \otimes X) \square (X \otimes X) +
(X \otimes Y) \square (Y \otimes X) = a \otimes a + b \otimes c = \Delta (a).$
\eproof

We can also recover algebraic relations of $Sl_q(2)$, except the $q$-determinant which is equal to zero,
i.e. $ad - q^{-1}bc =0$ instead of one. For checking the algebraic relations, we contract
the arrows $a,b,c,d$ into the vertices $\bar{a}=X^2, \bar{b}=XY, \bar{c}=YX, \bar{d}=Y^2$,
thanks to the usual product of $A$, (recall that such a product maps the arrow set
$A^{\otimes 2}$ into the vertex set $A$).
\begin{prop}
\label{important}
With $XY = \eta YX$, we obtain
$\bar{a}\bar{b}= \eta ^2 \bar{b}\bar{a}, \ \ \bar{c}\bar{b}=\bar{b}\bar{c}, \ \ \ \bar{a}\bar{c}= \eta ^2 \bar{c}\bar{a}, \ \
\bar{a}\bar{d}= \eta ^2 \bar{b}\bar{c}, \ \ \ \bar{b}\bar{d}= \eta ^2 \bar{d}\bar{b}, \ \ \ \bar{c}\bar{d}= \eta ^2 \bar{d}\bar{c},
\ \ \ \bar{a}\bar{d} - \bar{d}\bar{a} = (\eta ^2 - \eta ^{-2}) \bar{b}\bar{c}, \ \ \bar{a}\bar{d} -\eta ^2 \bar{b}\bar{c}=0.$
\end{prop}
\Proof
For instance, $\bar{a}\bar{b}= XX XY=  \eta ^2 XY XX = \eta ^2 \bar{b}\bar{a}$,
$\bar{c}\bar{d}= YX YY =  \eta ^2 YY YX = \eta ^2 \bar{d}\bar{c}$,
$\bar{a}\bar{d}= XX YY = \eta ^2 XYYX= \eta ^2 \bar{b}\bar{c}$, and so on. By setting $\eta^2 = q^{-1}$, we recover the usual algebraic relations
for $Sl_q(2)$, except the $q$-determinant which is equal to $0$.
\eproof
\section{The left part of $\F$}
The aim of this section is to prove that there exists a coassociative coalgebra $(\F, \tilde{\Delta})$ entangled to the usual coassociative coalgebra 
$(\F, \Delta)$, i.e. to prove that $(\F,\Delta, \tilde{\Delta})$ is a coassociative $L$-coalgebra.
By convention, we call $(\F, \tilde{\Delta})$, the left part of $(\F,\Delta, \tilde{\Delta})$. Its associated geometric support is:
\begin{center}
\includegraphics*[width=4cm]{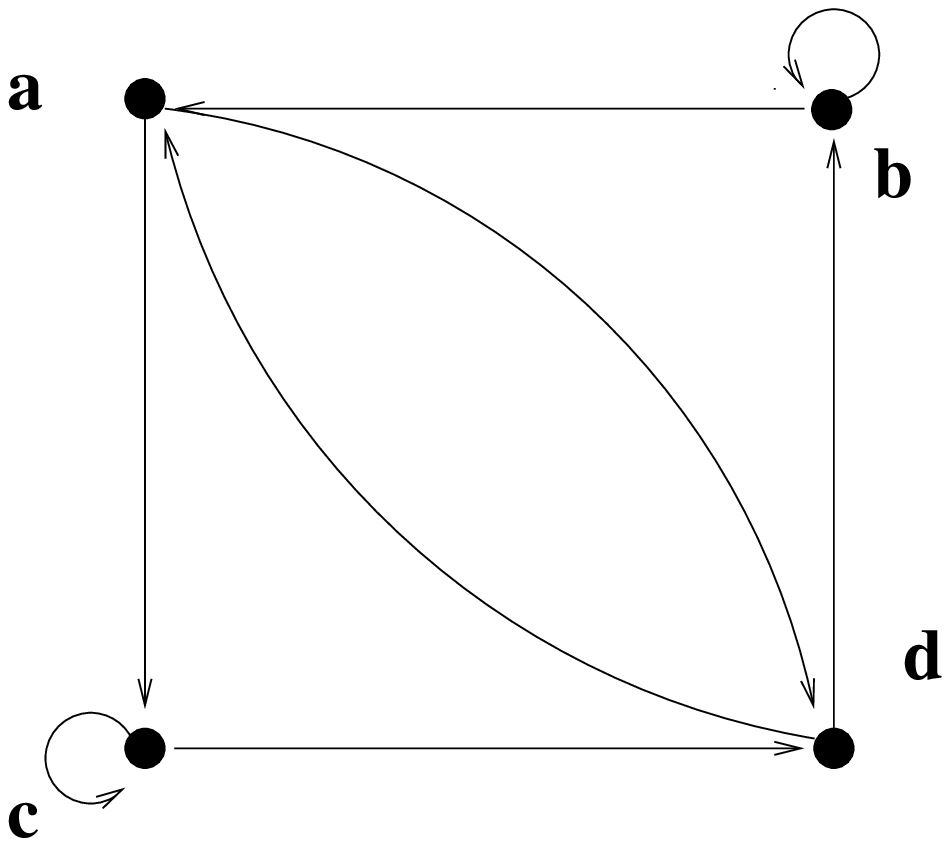}

\begin{scriptsize} \textbf{Geometric support associated with the left part $(\F, \tilde{\Delta})$.} \end{scriptsize}
\end{center}
This structure can be obtained, for instance, by inverting the map $\bra \cdot, \cdot \ket$ and the map $\per$ in the equations defining the usual
coproduct $\Delta$ of $\F$ obtained from the $(2,1)$-De Bruijn graph in proposition \ref{markov}.
The relations for the new coproduct $\tilde{\Delta}$ are,
$$\tilde{\Delta}b = b \otimes b + a \otimes d, \ 
\tilde{\Delta}c = c \otimes c + d \otimes a, \
\tilde{\Delta}a = b \otimes a + a \otimes c, \ 
\tilde{\Delta}d = c \otimes d + d \otimes b.$$
\Rk
In \cite{Coa}, we defined a matrix product $U \bar{\otimes} W$, where $U, W$ are two matrices and
$(U \bar{\otimes} W)_{ij} := \sum_k U_{ik} \otimes W_{kj}$. With:
\[
 U = \begin{pmatrix}
 a & b\\
 c& d
\end{pmatrix}, \ \
 \tilde{U} = \begin{pmatrix}
 b & a\\
 d& c
\end{pmatrix},
\]
the operation, $ \ \widetilde{} \ $, meaning here the permutation of the two columns, the definition of the two coproducts can be recovered.
Indeed,
$(U \bar{\otimes} U)_{ij} := \sum_k U_{ik} \otimes U_{kj}:= \Delta U_{ij}$ yields the right coproduct and the left coproduct
$\tilde{\Delta}$ can be recovered by computing $\widetilde{\tilde{U} \bar{\otimes} \tilde{U}}$.
\Rk
The linear map $\tilde{\epsilon}: \F \xrightarrow{} k$, such that:
$$\tilde{\epsilon}(a) = \tilde{\epsilon}(d) =0, \ \ \ \ \tilde{\epsilon}(b)=\tilde{\epsilon}(c) =1,$$
is a counit map for $\tilde{\Delta}$, i.e. $ (\tilde{\epsilon} \otimes id) \tilde{\Delta} =  (id \otimes \tilde{\epsilon}) \tilde{\Delta}= id.$
\begin{theo}
\label{theocoa}
The new coproduct $\tilde{\Delta}$ is coassociative and verifies the entanglement equation
$(\tilde{\Delta} \otimes id)\Delta = (id \otimes \Delta)\tilde{\Delta}$. Moreover the two coproducts verify also
$(\Delta \otimes id)\tilde{\Delta}= (id \otimes \tilde{\Delta}) \Delta$.
\end{theo}
\Proof
The coassociativity is straightforward. The proof of the two equations is also straightforward by using the matrix product defined above.
\eproof
\Rk
The axioms of coassociative co-dialgebras are not satisfied. 
\begin{defi}{[Chiral, Achiral]}
Let $(L,\Delta, \tilde{\Delta})$ be a coassociative $L$-coalgebra with right coproduct $\Delta$ and left coproduct $\tilde{\Delta}$.
The $k$-vector space $(L,\Delta, \tilde{\Delta})$ is {\it{chiral}} iff
$(\tilde{\Delta} \otimes id)\Delta = (id \otimes \Delta)\tilde{\Delta}$ and
$(\Delta \otimes id)\tilde{\Delta}\not= (id \otimes \tilde{\Delta}) \Delta$.
This means that the entanglement equation ``differentiates'' the left part $(L, \tilde{\Delta})$ from the right part $(L,\Delta)$. 
To the contrary, $(L, \Delta, \tilde{\Delta}) $ is said to be {\it{achiral}} iff
$(\tilde{\Delta} \otimes id)\Delta = (id \otimes \Delta)\tilde{\Delta}$ and
$(\Delta \otimes id)\tilde{\Delta}= (id \otimes \tilde{\Delta}) \Delta$, i.e. the axioms of an achiral coassociative $L$-coalgebra are globally invariant under the permutation $\tilde{\Delta}  \leftrightarrow \Delta$. More generally, an algebra
$(A, \ \bullet_1, \ldots, \bullet_n)$ equipped with $n$ operations $\bullet_1, \ldots, \bullet_n: A^{\otimes 2} \xrightarrow{} A$, verifying axioms $AX_1, \ldots AX_p$ is said to be achiral if $(A, \ \bullet_{\sigma(1)}, \ldots, \bullet_{\sigma(n)})$,
where $\sigma$ is a permutation, verifies also the same axioms, i.e. the axioms $AX_1, \ldots AX_p$ are globally invariant under the action of any permutations $\sigma$. The dualisation of this definition is straightforward.
\end{defi}
\Rk
In the theorem \ref{theocoa},
the entanglement equation is verified, even by inverting the r\^ole of the left and right coproducts. We say that $(\F,\Delta, \tilde{\Delta})$ is an
achiral $L$-coalgebra, since the left and right parts are entangled by the entanglement equation which do not differentiate
them. To the contrary, for instance, observe that the two coproducts $\delta_f$
and $\tilde{\delta}_f$ associated with an unital associative algebra, embed the algebra into a chiral $L$-bialgebra, see also \cite{dipt}.
\begin{theo}
The two geometric supports associated with the left  part $(\F,\tilde{\Delta})$ and the right part $(\F,\Delta)$ of the coassociative $L$-coalgebra $(\F,\Delta, \tilde{\Delta})$, obtained from the line-extension of the $(2,1)$-De Bruijn graph,
glued together yield the $(4,1)$-De Bruijn graph. Moreover the intersection of their arrow sets is empty.
\end{theo}
\Proof \textbf{[The gluing of the left and right parts of $(\F,\Delta, \tilde{\Delta})$]}
When we glue the geometric support associated with $(\F,\Delta)$ with its left part $Gr(\F, \tilde{\Delta})$, we obtain the $(4,1)$-De Bruijn directed graph.
\begin{center}
\includegraphics*[width=3.5cm]{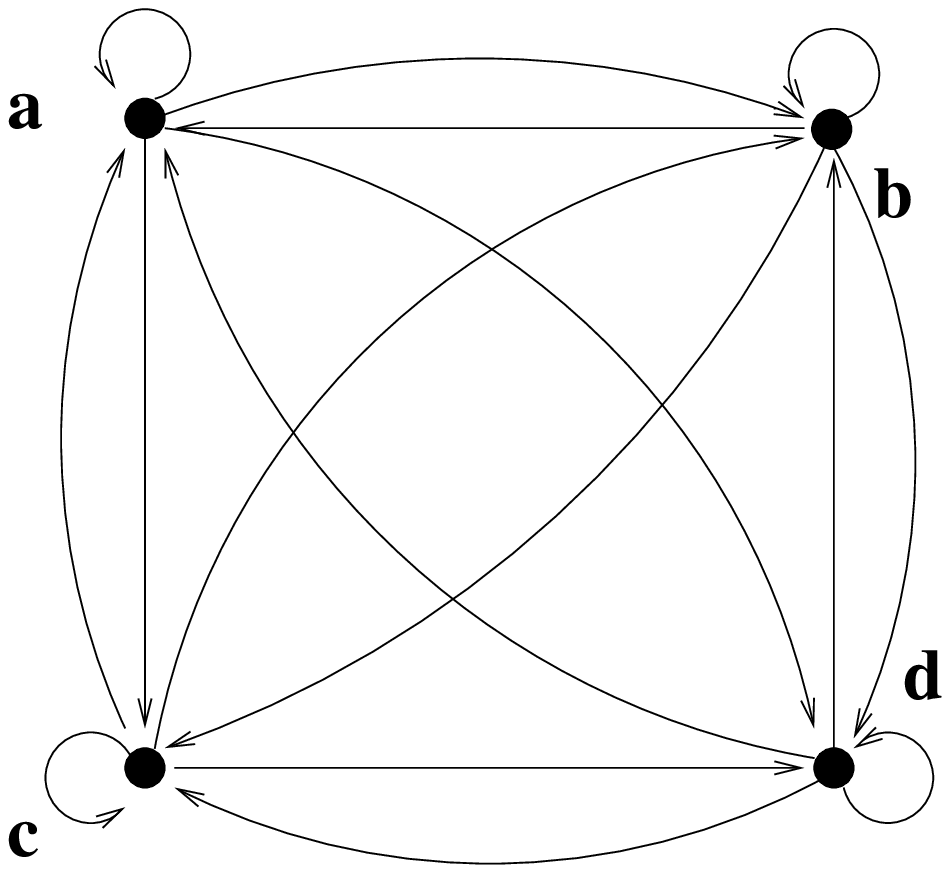}

\textbf{\begin{scriptsize} The $(4,1)$-De Bruijn graph. \end{scriptsize}}
\end{center}
The intersection of the arrow sets of these two graphs is empty.
\eproof

\noindent
The $(4,1)$-De Bruijn graph is tiled with two (geometric supports of) entangled coassociative coalgebras, the entanglement being achiral. \subsection{Consequences}
\subsubsection{An example of achiral $L$-Hopf algebra}
\begin{prop}
If $(\F,\Delta, \tilde{\Delta})$ is also a unital algebra and $a,b,c,d$ commute pairewise. Then $(\F,\Delta, \tilde{\Delta})$ will be an achiral $L$-bialgebra.
Moreover, if $ad - bc =1$,
define the linear map $S_\F$ which maps $a \mapsto d$,
$d \mapsto a$, $b \mapsto -b$ and $c \mapsto -c$ and the linear map $\tilde{S}_\F$ which maps
$b \mapsto -c$, $c \mapsto -b$, $a \mapsto a$, $d \mapsto d$, then the right (resp. left) part of $(\F,\Delta, \tilde{\Delta})$ is a Hopf algebra, i.e
$(\F,\Delta, \tilde{\Delta})$ is an achiral $L$-Hopf algebra.
\end{prop}
\Proof
Straightforward.
\eproof
\subsubsection{An example of chiral $L$-bialgebra}
Consider the algebra generated by $x,p$ and $g$
such that $pg=gp$, $xg = \eta gx$, $xp = \mu px$, with  $\eta, \mu $ two scalars different from zero.
Define the (left) part by: $\tilde{\Delta} x = x \otimes p + g \otimes x$,  $\tilde{\Delta} p = p \otimes p$, $\tilde{\Delta} g = g \otimes  g$
with counit $\tilde{\epsilon}$: $x \mapsto 0$ and $p,g \mapsto 1$. Define
the (right) part by: $\Delta x = x \otimes x$,  $\Delta p = p \otimes x$, $\Delta g = g \otimes  x$.
\begin{prop}
The two coproducts are coassociative and verify the entanglement equation $(id \otimes \Delta) \tilde{ \Delta} =
(\tilde{ \Delta} \otimes id) \Delta $. Moreover the two coproducts $\Delta$, $\tilde{\Delta}$ and the counit $\tilde{\epsilon}$ are homomorphisms.
\end{prop}
\Proof
Tedious but straightforward by noticing that it is a consequence of the previous subsection on the left part
of $\F$, by setting formally $b=0$. For instance, we check the entanglement equation on $x$. We get,
$x \xrightarrow{\Delta} x \otimes x \xrightarrow{\tilde{\Delta} \otimes id}x \otimes p  \otimes x +  g \otimes x \otimes x$ and
$x \xrightarrow{\tilde{\Delta}} x \otimes p + g  \otimes x  \xrightarrow{id \otimes \Delta} x \otimes p \otimes x + g  \otimes x \otimes x$.
For the homomorphism property, we get for instance  $\Delta xg =  \Delta x \cdot \Delta g = xg \otimes xx$ and
$\Delta gx =  \Delta g \cdot \Delta x = gx \otimes xx$, thus $\Delta xg = \eta \Delta gx$ and so forth.
\eproof

\subsubsection{Consequences for $Sl_q(2)$}
\begin{prop}
With the usual algebraic relations of the Hopf algebra $Sl_q(2)$, whose usual product will be denoted by $m$, the linear map
$\tilde{S}: Sl_q(2) \xrightarrow{}Sl_q(2)$ defined by $b \mapsto -q^{-1}c; c \mapsto -qb; a \mapsto a; d \mapsto d$,
verifies $m(id \otimes \tilde{S}) \tilde{\Delta} =m (\tilde{S} \otimes id) \tilde{\Delta} = 1 \cdot \tilde{\epsilon}.$
Moreover, the map $\tilde{S}$ can be extended to an unital algebra map. As $1,a,d$ and $bc$ are fixed points of $\tilde{S}$,
the $q$-determinant of the matrix $U$ is invariant by applying $\tilde{S}$.
\end{prop}
\Proof
Let us check the antipode property. For that, multiply $\tilde{U}$ by the following matrix
\[
 \tilde{S}(\tilde{U}) = \begin{pmatrix}
 -q^{-1}c & a\\
 d& -qb
\end{pmatrix}, \]
which is the matrix obtained from $\tilde{U}$ by the left antipode $\tilde{S}$. We compute
\[ \widetilde{\tilde{U} \tilde{S}(\tilde{U})} =\begin{pmatrix}
 ba -qab & ad -q^{-1}bc\\
 da -qcb & cd -q^{-1}dc
\end{pmatrix}. \]
This matrix must be equal to:
\[ \tilde{\epsilon}(U)=\begin{pmatrix}
  \tilde{\epsilon}(a) =0 & \tilde{\epsilon}(b) =1 \\
 \tilde{\epsilon}(c) =1 & \tilde{\epsilon}(d) =0
\end{pmatrix}, \]
which is the case since we know that $da = 1 + qbc$. We obtain the same result if we consider $\widetilde{\tilde{S}(\tilde{U})\tilde{U}}$.
The extension of the definition of $\tilde{S}$ to a homomorphism is straightforward.
\eproof
\Rk
The algebraic relations of $Sl_q(2)$ do not embed its left part into a bialgebra. This fact will be a motivation for introducing coassociative manifolds in \cite{dipt}.
\subsubsection{The left part of $SU_q(2)$}
A consequence of what was done for the coassociative coalgebra $\F$ is the existence of a left part for the Hopf algebra $SU_q(2)$.
We recall that $SU_q(2)$ is an $*$-algebra with two generators, $a$ and $c$ such that
$ac = qca, \ ac^* = qc^*a, \ cc^* = c^*c, \ a^*a + c^*c =1, \ aa^* + q^2 cc^* =1,$ with $q$ a real different from $0$.
The coproduct is given by $\Delta_1 a = a \otimes a - qc^* \otimes c$ and $\Delta_1 c = c \otimes a + a^* \otimes c$ and
the counit is $\epsilon_1(a) =1$ and $\epsilon_1(c) = 0$.
\begin{prop}
There exists a coassociative coalgebra which is achiral entangled to
the usual coassociative coalgebra part of $SU_q(2)$.
\end{prop}
\Proof
Define the (left) coproduct $\tilde{\Delta}_1$
by the $*$-linear map $a \mapsto \tilde{\Delta}_1 a := c^* \otimes a + a \otimes c$ and
$c \mapsto \tilde{\Delta}_1 c := c \otimes c - q^{-1} a^* \otimes a$.
This coproduct is coassociative. Moreover, we have
$(\tilde{\Delta}_1 \otimes id) \Delta_1 = (id \otimes \Delta_1) \tilde{\Delta}_1$ and
$(\Delta_1 \otimes id) \tilde{\Delta}_1 = (id \otimes \tilde{\Delta}_1) \Delta_1.$
For instance, $(\Delta_1 \otimes id) \tilde{\Delta}_1 a = c^* \otimes a^* \otimes a + a \otimes c^* \otimes  a + a \otimes a \otimes c
- qc^* \otimes c \otimes c$ and
$(id \otimes \tilde{\Delta}_1) \Delta_1 a = a \otimes c^* \otimes a
+ a \otimes a \otimes c -q(c^* \otimes c \otimes c  - q^{-1}c^* \otimes a^* \otimes a),$ and so forth.
The $*$-linear map $\tilde{\epsilon}_1$ is defined by $c \mapsto 1$ and $a \mapsto 0$
is the (left) counit.
\eproof
\Rk
The usual algebraic relations of $SU_q(2)$ do not embed the new coproduct and the new counit into homomorphisms.

\subsection{Splitting the coassociative coproducts and achiral $L$-coalgebras}
The example of dendriform algebras \footnote{See the introduction or \cite{Loday}.} shows that there
exists associative algebras whose associative product can split into two (a priori non associative) products. In the case of achiral
$L$-coalgebras, we have also a splitting of a coassociative coproduct into two coassociative ones.
Let us see how it works.

\noindent
Recall that an achiral $L$-coalgebra $(L,\Delta,\tilde{\Delta})$ has two coproducts $\Delta$ and $\tilde{\Delta}$ which verify the axioms:
\begin{enumerate}
\item {$\Delta$ and $\tilde{\Delta}$ are coassociative.}
\item{$(\tilde{\Delta} \otimes id) \Delta = (id \otimes \Delta) \tilde{\Delta}.$}
\item{$(\Delta \otimes id) \tilde{\Delta} = (id \otimes \tilde{\Delta}) \Delta,$}
\end{enumerate}
i.e. the axioms are globally invariant by the permutation $\Delta \leftrightarrow \tilde{\Delta}$.
\begin{prop}
Suppose $L$ is an achiral $L$-coalgebra.
Let $u,v,w,z \in k$.
These axioms are invariant under the transformation
$ \begin{pmatrix}
\Delta \\
\tilde{\Delta}
\end{pmatrix}
\mapsto
\begin{pmatrix}
u & v \\
w & z
\end{pmatrix}
\begin{pmatrix}
\Delta \\
\tilde{\Delta}
\end{pmatrix}$. For instance, the coproduct $\Delta ' := u \Delta + v \tilde{\Delta}$
is still coassociative.
\end{prop}
\Proof
Straightforward. 
\eproof
\section{Tiling the $(n^2,1)$-De Bruijn co-dialgebra with $n$ coassociative coalgebras}
We can generalize the previous procedure to any coassociative coalgebras whose geometric supports are obtained by line-extension of the $(n,1)$-De Bruijn graphs
viewed as geometric supports of Markov co-dialgebras.
Fix $n>1$, the number of vertices of the $(n,1)$-De Bruijn graph, denoted by $U_i$, $i=1, \ldots n$.
We have seen in proposition \ref{tyu} that the line-extension of such a graph yields a geometric support of a coassociative coalgebra whose coproduct is denoted by $\Delta$. The coproduct was recovered
by computing $(U \bar{\otimes} U)$, where $U_{ij}$ is the vertex associated with the arrow going from the vertex $U_i$ to $U_j$ in the $(n,1)$-De Bruijn
graph. Denote by $p$, the shift which maps $j $ into $j+1 \mod n$, for all $j = 1, \ldots, n$ .
Let $\alpha$ be an integer equal to $0, \ldots, n-1$. We denote by $\mathcal{P}{^\alpha}(U)$,
(resp. $\mathcal{P}^{- \alpha}(U)$),
the matrix obtained from $U$ by letting the shift $p^\alpha$, (resp. $p^{-\alpha}$) acts on the columns of $U$, i.e. $\mathcal{P}^\alpha
(U)_{ij} := U_{ip^\alpha(j)}$
and $\mathcal{P}^{-\alpha}(U)_{ij} := U_{ip^{-\alpha}(j)}$.
\begin{lemm}
\label{lmm}
Let $A,B$ be two $n$ by $n$ matrices and $\alpha$ be an integer equal to $0, \ldots, n-1$. We get $\mathcal{P}^\alpha(AB) = A\mathcal{P}^\alpha(B)$ and  $\mathcal{P}^{- \alpha}(AB)
= A \mathcal{P}^{- \alpha}(B) $.
Therefore,
we obtain, $\mathcal{P}^\alpha (U)U = \mathcal{P}^\alpha (\mathcal{P}^\alpha (U) \mathcal{P}^{- \alpha }(U))$.
\end{lemm}
\Proof
Notice that $\mathcal{P}^\alpha(AB)_{ij}=(AB)_{ip^\alpha(j)}= \sum_k \ A_{ik}B_{kp^\alpha(j)}= (A\mathcal{P}^\alpha(B))_{ij}$,
which proved the first equality. The sequel is now straightforward.
\eproof
\Rk
This lemma is also valid by replacing the usual product by  $\bar{\otimes}$.
\begin{defi}{[Coproducts]}
As $\mathcal{P}^0 = id$, we rename the usual coproduct $\Delta$ by $\Delta_{[0]}$. Its explicit definition,
as we have seen, is closely related to the matrix $U$.  We define
also the coproducts $\Delta_{[\alpha]}$, for all $\alpha$, by
$\Delta_{[\alpha]}U_{ij} =(\mathcal{P}^\alpha (\mathcal{P}^\alpha (U) \bar{\otimes}\mathcal{P}^{- \alpha }(U)) _{ij}
=(\mathcal{P}^\alpha (U)U)_{ij}$. Denote by $(\F_n, \Delta_{[\alpha]})$, $0 \leq \alpha \leq n-1$, the coalgebra so obtained. (Observe that $\F_2 := \F$.)
\end{defi}
\begin{theo}
The new coproducts $\Delta_{[\alpha]}$ are coassociative. Moreover, for all $\alpha, \beta = 0, \ldots, n-1$,  $\Delta_{[\alpha]}$ and
$\Delta_{[\beta]}$ obey the entanglement equation.
\end{theo}
\Proof
We define for two matrices $A,B$ the following product $A*_\alpha B =(\mathcal{P}^\alpha (\mathcal{P}^\alpha (A) \bar{\otimes}\mathcal{P}^{- \alpha }(B)) $. From the straightforward equality $A*_\alpha (B*_\alpha C) = (A*_\alpha B)*_\alpha C$, we obtain
in the case where $A=B=C=U$, the coassociativity equation
$(\Delta_{[\alpha]} \otimes id)\Delta_{[\alpha]}U_{ij}=(id \otimes \Delta_{[\alpha]}) \Delta_{[\alpha]}U_{ij}$.
To prove that two coproducts obey the entanglement equation, we have to show that
$(\Delta _{[\alpha]}\otimes id)\Delta_{[\beta]} U_{ij}=(id \otimes \Delta_{[\beta]}) \Delta_{[\alpha]}U_{ij}$, that is
$U *_\alpha (U *_\beta U) = (U *{_\alpha} U) *_\beta U$, which is also straightforward.
\eproof
\Rk
We have showed that in the case of a coassociative coalgebra obtained by line-extension of the $(n,1)$-De Bruijn graphs,
we can construct others coassociative coalgebras whose coproducts verify the entanglement equation. Precisely,
the coassociative coalgebras $(\F_n, \Delta_{[\alpha]})$ and $(\F_n, \Delta_{[\beta]})$, $0 \leq \alpha,\beta \leq n-1$, are entangled, the entanglement being achiral. Therefore, the $k$-vector space $(\F_n, \Delta_{[0]}, \ldots, \Delta_{[n-1]})$ is achiral\footnote{In \cite{dipt}, such algebraic objects will be called coassociative manifolds.}. 
\begin{prop}
Let $(C, \Delta_{[0]}, \ldots, \Delta_{[n-1]})$ be a $k$-vector space equipped with $n$ linear maps $\Delta_{[i]}: C \xrightarrow{} C^{\otimes 2}$, $i=0, \ldots n-1$, such that:
$$ (\Delta_i \otimes id) \Delta_j = (id \otimes \Delta_j ) \Delta_i, \ \ i,j=0, \ldots, n-1.$$
For $i=0, \ldots n-1$, set $x_i$ the vector equal to $^t(0,\ldots, 0,\Delta_{i},0, \ldots, 0)$ and fix $Z \in M_n(k)$, a $n$
by $n$ matrix. Then $(C, \Delta'_{[0]}:=Zx_0, \ldots, \Delta'_{[n-1]}:=Zx_{n-1})$ verify also:
$$ (\Delta'_i \otimes id) \Delta'_j = (id \otimes \Delta'_j ) \Delta'_i, \ \ i,j=0, \ldots, n-1.$$
\end{prop}
\Proof
Straightforward.
\eproof
\Rk
Observe that this theorem allow us to produce splittings of associative laws into associative ones.
\Rk
Applying this proposition to the tiling of the $(3,1)$-De Bruijn graph gives examples of cubical cotrialgebras, i.e. a coalgebra equipped with 3 coproducts $\Delta_i$, $i=0,1,2,$ such that:
$$ (\Delta_i \otimes id) \Delta_j = (id \otimes \Delta_j ) \Delta_i, \ \ i,j=0,1,2.$$
The notion of cubical trialgebra is defined in \cite{LodayRonco}. We recover easily cubical trialgebras from cubical cotrialgebras by considering the convolution products. The operad $Tricub$ on one generator associated with cubical trialgebra is Koszul and self-dual. The operad associated with the so-called
hypercube $n$-algebra, i.e. a $k$-vector space
equipped with $n$ products verifying:
$$ (x \bullet_i y) \bullet_j z = x \bullet_i (y \bullet_j z), \ \ x,y,z \in A, \ i,j=0, \ldots, n-1,$$
is conjectured to be Kozul and self-dual (observe
that there are $n^2$ operations and two possible choices of parentheses, thus $2n^2 - n^2 =n^2$).
\Rk
The counit $\epsilon_{[\alpha]}$, associated with the coassociative coproduct $\Delta_{[\alpha]}$, is obviously defined by
$U_{ip^\alpha(i)} \mapsto 1$ and $U_{ip^\alpha(j)} \mapsto 0$ if $i \not= j$.
\begin{theo}
The intersection of the arrow sets of the geometric supports of $(\F_n, \Delta_{[\alpha]})$, $0 \leq \alpha \leq n-1$, are empty, i.e. $Gr((\F_n, \Delta_{[\alpha]}))_1 \cap Gr((\F_n, \Delta_{[\beta]}))_1 = \emptyset$, $0 \leq \alpha, \ \beta \leq n-1$ and $\alpha \not= \beta$. Gluing them yields the
$(n^2,1)$-De Bruijn directed graph, i.e. $\bigcup_{0 \leq \alpha \leq n-1} Gr((\F_n, \Delta_{[\alpha]})) = D_{(n^2,1)}$.
\end{theo}
\Proof
We call $(U_{ij})_{i,j = 1, \ldots, n}$, the vertices of the $(n^2,1)$-De Bruijn graph.
Let us prove that the gluing of all the geometric supports of $(\F_n, \Delta_{[\alpha]})$, $0 \leq \alpha \leq n-1$, yields the
$(n^2,1)$-De Bruijn directed graph. Every arrow of the $(n^2,1)$-De Bruijn directed graph can be described by
$U_{ik} \otimes U_{lj}$, with $k,l,i,j = 1, \ldots, n$.
As the shift is one-to-one, there exists an unique integer $\alpha$
such that $p^\alpha(l) = k$ i.e. $U_{ik} \otimes U_{lj} = U_{ip^\alpha(l)} \otimes U_{lj}$,
i.e this arrow belongs to the definition of the coproduct $\Delta_[\alpha]$, see lemma \ref{lmm}. Therefore, the $(n^2,1)$-De Bruijn graph is included into the gluing of
the geometric support of the $n$ coassociative coalgebras. The reversal is obvious.

\noindent
Fix $\alpha$ and $\beta$
two different integers. As the shift $p^{\beta - \alpha}$ has no fixed point,
no arrow defined in the coproduct $\Delta_{[\alpha]}$ is present in the definition of the coproduct of
$\Delta_{[\beta]}$. Therefore, the intersection of arrow sets of $Gr(\F_n, \Delta_{[\alpha]})$ and $Gr(\F_n, \Delta_{[\beta]})$ is empty.
\eproof
\Rk
Via their geometric supports, we get:\\
$(n,1)$-De Bruijn co-dialgebra $ \xrightarrow{\textsf{Line-extension}}$ $n$ coassociative coalgebras \footnote{
entangled by the achiral entanglement equation.}
$\xrightarrow{\textsf{Gluing}}$ $(n^2,1)$-De Bruijn co-dialgebra.

\noindent
We have yielded a tiling of the $(n^2,1)$-De Bruijn graph into $n$ coassociative coalgebras, each 
$(\F_n, \Delta_{[\alpha]}, \Delta_{[\beta]})$, $0 \leq \alpha,\beta \leq n-1$, being an achiral coassociative $L$-coalgebra. The case $n=1$ is trivial. Indeed, the $(1,1)$-De Bruijn graph and
its line-extension are loops and a loop is coassociative, since it is of the form $x \mapsto x \otimes x$.
\NB
Fix $n \geq 1$.
The link between the Markovian coproducts of the coassociative codialgebra $(k(D_{(n,1)})_0, \Delta_M, \tilde{\Delta}_M)$, associated with the $(n,1)$-De Bruijn graph and the coassociative coproduct of $(\F_n, \Delta_{[0]})$ can be seen as follows.
Recall that $D_{(n,1)}$ is the $(n,1)$-De Bruijn graph, with vertex set $(D_{(n,1)})_0:=\{U_i; \ 1 \leq i \leq n \}$. Consider its free $k$-vector space $k(D_{(n,1)})_0$. Its Markovian coproducts verify $\Delta_M U_i:= U_i \otimes \sum_k U_k$ and 
$\tilde{\Delta}_M U_i:= \sum_k U_k \otimes  U_i$.
Set for all $i,j=0, \ldots, n-1$, $U_{ij}:= U_i \otimes U_j$. Define the map $\boxminus: k(D_{(n,1)})_0^{\otimes 2}\times k(D_{(n,1)})_0^{\otimes 2} \xrightarrow{} k(D_{(n,1)})_0^{\otimes 2} \otimes k(D_{(n,1)})_0^{\otimes 2}$ such that
$(z_1 \otimes z_2) \boxminus (z_3 \otimes z_4):=
(id \otimes \Psi \otimes id)(z_1 \otimes z_2 \otimes z_3 \otimes z_4)$, where $\Psi(z_i \otimes z_j)= z_i \otimes z_i$ if $i=j$, and 0 otherwise.
Observe that:
$$ \Delta_M(U_i) \boxminus \tilde{\Delta}_M(U_j) := \Delta_{[0]} U_{ij}, \ \ i,j=0, \ldots n-1. $$

\noindent
As in $L$-coalgebra theory, the notion of algebraic product is not assumed,  we end this section by a propositon
on the Leibniz coderivation and an example applied on $\F$.

\begin{defi}{[Coderivation]}
Let $C$ be a coassociative coalgebra with coproduct $\Delta$.
A linear map $D: \ C \xrightarrow{} C $ is called a (Leibniz) {\it{coderivation}} with respect to the coproduct $\Delta$ iff it verifies:
$$ \Delta D = (id \otimes D) \Delta + (D \otimes id) \Delta.$$
\end{defi}
\begin{prop}
Let $C$ be a coassociative coalgebra with coproduct $\Delta$, such that for all the elements $U_{ij} \in C$,
$\Delta U_{ij} = \sum_{k=1}^n U_{ik} \otimes U_{kj}$, where $n$ is a fixed integer.
Define the linear map $D$ by $U_{ij} \mapsto \sum_k ^n U_{kj} - U_{ik}$. Then $D$ is a Leibniz coderivation.
\end{prop}
\Proof
Fix $U_{ij}$. We have  $(id \otimes D) \Delta + (D \otimes id) \Delta(U_{ij} ) = \sum_{k,l =1}^n U_{lk} \otimes U_{kj} - U_{ik} \otimes U_{kl}$
which is equal to $\Delta (\sum_{l=1}^n U_{lj} - U_{il}).$
\eproof
\begin{exam}{[A coderivation for $(\F,\Delta, \tilde{\Delta})$]}\\
For $(\F,\Delta)$, we get $D(a) = c-b = -D(d), \ D(b)= d-a = -D(c)$. As $(\F,\Delta, \tilde{\Delta})$ is a $L$-coalgebra, we yield also a coderivation $\tilde{D}$ for
$(\F, \tilde{\Delta})$.  Define $\tilde{D} = D$. A straightforward computation shows
that $D$ is also a coderivation with respect to the coproduct $\tilde{\Delta}$.
\end{exam}

\section{Conclusion}
The first main result obtained in this paper is the possibility to construct from directed graphs, families of $L$-cocommutative coassociative co-dialgebras and
therefore, via convolution products families of associative dialgebras. The second main result is the possibility to
recover from the line-extension of the geometric supports of these coassociative co-dialgebras, known coassociative coalgebras and in the case
of the $(n^2,1)$-De Bruijn graphs, $n > 0$, to obtain a tiling of these Markovian objects by $n$ (geometric supports of) coassociative coalgebras.
An important notion, called the achirality has been put forward. We have shown that actions of $M_n(k)$ on the coproducts defining the tiling of the $(n^2,1)$-De Bruijn graph let globally invariant the relations between them.
We gave consequences of such tilings and found
examples of cubical trialgebras and more generally, examples of hypercube $n$-algebras. In addition, this allowed us to construct associative laws which split into several associative ones.

This paper has been pursued in \cite{dipt}. In 
\cite{dipt}, the notion of codipterous coalgebras and 
pre-dendriform coalgebras \footnote{These notions have been discovered by J-L Loday and M. Ronco \cite{LR1} and rediscovered independently, via graph theory by the author \cite{LL}.} are established. These spaces
constructed from coassociative coalgebra theory 
extend the notions developed so far and are the elementary boxes of coassociative manifolds. Via 
these notions, we construct Poisson algebras,
dendriform algebras, associative dialgebras (which are not Markovian), associative
trialgebras \cite{LodayRonco}. Notably,
the tilings constructed so far will yield examples
of coassociative manifolds \cite{dipt}.

\noindent
\textbf{Acknowledgments:}
The author wishes to thank Dimitri Petritis for useful discussions and fruitful advice for the
redaction of this paper and to S. Severini for pointing him the precise definition of the De Bruijn graphs.\\

\bibliographystyle{plain}
\bibliography{These}

\end{document}